\documentclass[reqno]{amsart}

\usepackage{amsfonts,amsmath,amsthm,amssymb, amscd}
\usepackage{enumerate, url, desclist}
\usepackage[english]{babel}
\usepackage[latin1]{inputenc}

\newtheorem{thm}{Theorem}[section]
\newtheorem*{thm*}{Theorem}
\newtheorem{cor}[thm]{Corollary}
\newtheorem{prop}[thm]{Proposition}
\newtheorem{lemma}[thm]{Lemma}
\newtheorem{step}{Step}
\theoremstyle{definition}
\newtheorem{defin}{Definition}[section]
\newtheorem*{rem}{Remark}
\newtheorem*{sketch}{Sketch of strategy}
\newtheorem*{ack}{Acknowledgement}

\newenvironment{mainthm}[1]{\renewcommand{\thethm}{#1}\begin{thm}}{\end{thm}}

\numberwithin{equation}{section}

\newcommand{\N}{\mathbb{N}}
\newcommand{\Z}{\mathbb{Z}}
\newcommand{\Q}{\mathbb{Q}}
\newcommand{\F}{\mathbb{F}}

\newcommand{\R}{\mathcal{R}}
\renewcommand{\P}{\mathcal{P}}
\newcommand{\gm}{\mathbb{G}_m}

\newcommand{\ktl}{\widetilde{k}}
\newcommand{\ve}[1]{\mathbf{#1}}
\newcommand{\cl}[1]{\overline{#1}}
\newcommand{\mc}[1]{\mathcal{#1}}

\newcommand{\h}[1]{\mc{H}(#1)}
\newcommand{\res}[1]{\text{\rule[-1ex]{0.1ex}{2.8ex}}_{\, #1}}

\DeclareSymbolFont{symbolsC}{U}{pxsyc}{m}{n} \SetSymbolFont{symbolsC}{bold}{U}{pxsyc}{bx}{n}
\DeclareFontSubstitution{U}{pxsyc}{m}{n} \DeclareMathSymbol{\nequiv}{\mathrel}{symbolsC}{"2E}

\DeclareMathOperator{\Gal}{Gal} \DeclareMathOperator{\Orb}{Orb}  \DeclareMathOperator{\id}{I}

\title{Equations in the Hadamard ring of rational functions}
\author{Andrea Ferretti}
\address{Andrea Ferretti
\newline
Dipartimento di matematica - Università "La Sapienza"
\newline
Piazzale Aldo Moro, 2 - 00185 Roma, Italy}
\email{ferretti@mat.uniroma1.it}
\author{Umberto Zannier}
\address{Umberto Zannier
\newline
Scuola Normale Superiore
\newline
Piazza dei Cavalieri, 7 - 56126 Pisa, Italy}
\email{u.zannier@sns.it}
\date{January 26, 2007}

\begin{document}

\begin{abstract}
Let $k$ be a number field. It is well known that the set of sequences composed by Taylor coefficients of rational functions over $k$ is closed under component-wise operations, and so it can be equipped with a ring structure. A conjecture due to Pisot asks if (after enlarging the field) one can take $d^{th}$ roots in this ring, provided $d^{th}$ roots of coefficients can be taken in $k$. This was proved true in a preceding paper of the second author; in this article we generalize this result to more general equations, monic in $Y$, where the former case can be recovered for $g(X,Y)=X^d-Y=0$. Combining this with the \emph{Hadamard quotient theorem} by Pourchet and Van der Poorten, we are able to get rid of the monic restriction, and have a theorem that generalizes both results.
\end{abstract}

\maketitle

\section{Introduction}

Let $k$ be a field of characteristic $0$. We define a \emph{recurrence sequence} to be a sequence $\{ a(n) \}_{n \in \N} \subset \cl{k}$ satisfying
\begin{equation*}
a(n+m)+c_{m-1} a(n+m-1) +\dots +c_0 a(n) = 0
\end{equation*}
for each $n \geq 0$, for some fixed $c_0, \dots c_{m-1} \in k$. When $m$ is minimal, the polynomial
\begin{equation*}
q(T)=T^m + c_{m-1} T^{m-1} +\dots +c_0
\end{equation*}
is said to be \emph{associated} with the recurrence, and its roots $\alpha_i$ are by definition the \emph{roots} of the recurrence.

On the other hand, whenever a rational function $R \in k(x)$ is defined in $0$, Taylor coefficients may be taken, setting as usual $s_k=R^{(k)}(0)/k!$. It is well known that in the ring of formal power series the equality $R(x)=\sum s_k x^k$ holds, and it is easy to show that a sequence $\{ s_k \}$ represents a rational function if, and only if, it satisfies a linear recurrence except for a finite number of terms. In this case we call it a \emph{rational power series}.

Now, it is well known (for all these facts see \cite{BetterKnown}) that recurrence sequences are characterized by an explicit closed form, given by \emph{exponential polynomials}
\begin{equation*}
a(n)=\sum_{i=1}^{m} A_i(n) \alpha_i^n
\end{equation*}
where $A_i \in \cl{k}[x]$ and $\alpha_i \in \cl{k}$ (the $\alpha_i$ are in fact the roots of the recurrence). Since the sum and products of exponential polynomials are sequences of the same kind, it follows that the set of recurrence sequences (or equivalently the set of rational power series) is closed under component-wise sum and product. This leads us to make the following
\begin{defin}
The \emph{Hadamard ring} over the field $k$ is the set of formal power series with coefficients in $k$ which represent a rational function, equipped with component-wise operations. Equivalently it can be thought as the set of sequences from $k$ definitively satisfying a linear recurrence. It is denoted by $\h{k}$. Whenever $a \in \h{k}$ we denote by $a(n)$ its $n$-th coefficient (or its $n$-th term, if you think of recurrence sequences).
\end{defin}

Suppose now that we want to solve an algebraic equation in $\h{k}$: the first attempt is to solve it in the larger set $k[[x]]$, which we regard as a ring under component-wise sum and product of coefficients. This, in turn, amounts to solve infinitely many equations in the field $k$.

The case we are interested in is when $k$ is a number field, and we shall assume this from now on. We shall also identify a formal power series with the sequence of its coefficients. With this terminology Zannier proves the following theorem, solving a conjecture of Pisot:
\begin{thm*}[\cite{Z00}]
Let $k$ be a number field and let $\sum b(n) x^n \in \h{k}$. Suppose that for all $n$ the equation $Y^d=b(n)$ has a solution in $k$. Then there exists a finite extension $k'/k$ such that the same equation has a solution in $\h{k'}$. In other words we may choose $d$-th roots for the $b(n)$ so that they satisfy themselves a linear recurrence.
\end{thm*}
Another classical result for the problem of solving equations in this ring is the \emph{Hadamard quotient theorem} (proved in \cite{Pour} and \cite{VP88}, but see also \cite{Ru86} for a detailed account), which deals with linear equations.
\begin{thm*}[Hadamard quotient]
Let $F$ be a field of characteristic zero and let $b(n), c(n) \in \h{F}$. Let  $(a_n)$ be a sequence whose elements are in a subring $R$ of $F$ which is finitely generated over $\Z$, and suppose that $a_n=b(n)/c(n)$ whenever the quotient is defined. Then there exists an element $a(n) \in \h{F}$ such that $a(n)=a_n$ for every $n$ such that $c(n) \neq 0$.
\end{thm*}

In this paper we generalize these results, providing a solution to a more general conjecture of Van der Poorten (\cite{Rocky}).
\begin{thm} \label{main}
Let $k$ be a number field, $b_0, \dots, b_{d-1} \in \h{k}$, and consider the equation
\begin{equation} \label{problem}
Y^d+b_{d-1}(n)Y^{d-1}+ \dots +b_0(n)=0.
\end{equation}
Suppose \eqref{problem} has a solution for all $n$; then there exists a finite extension $k'/k$ such that the same equation has a solution in $\h{k'}$.
\end{thm}
It will be convenient to restate theorem \ref{main} in terms of exponential polynomials; moreover we may assume that the $b_j(n)$ have roots contained in the same finite set $\{ \beta_1, \dots, \beta_m \}$.
\begin{mainthm}{\ref{main}}[$2^{nd}$ form] \label{main2}
Let $k$ be a number field and for $j=0, \dots, d-1$ let
\begin{equation*}
b_j(n)=\sum_{i=1}^{m} B_{i,j}(n) \beta_i^n
\end{equation*}
be exponential polynomials, with $B_{i,j} \in k[x]$ and $\beta_i \in k$ for all $i,j$. Suppose that for every $n$ the equation
\begin{equation}
Y^d+b_{d-1}(n)Y^{d-1}+ \dots +b_0(n)=0 \tag{\ref{problem}}
\end{equation}
has a solution $a_n \in k$. Then there exists an exponential polynomial $a(n)$ with coefficients in a finite extension of $k$ that satisfies \eqref{problem} identically.
\end{mainthm}
\addtocounter{thm}{-1}

\begin{rem}
Of course one can relax the hypothesis requiring that the equations have solution in a fixed finite extension of $k$. Actually we will enlarge $k$ in the course of the proof without further comment.
\end{rem}
\begin{rem}
One can use the techniques of reduction of Rumely and Van der Poorten (\cite{RuVdP}, \cite{Ru86}) to deduce from theorem \ref{main} an analogous statement for a field $k$ finitely generated over $\Q$. We omit this verification, which is substantially straightforward after the quoted papers. See also the paper of Corvaja \cite{Corvaja} for a somewhat different deduction.
\end{rem}

Since our proof will involve an induction, it will be convenient to state and prove the following stronger form of theorem \ref{main2}.
\begin{thm}\label{main3}
Suppose that for each arithmetic progression $\mathfrak{A}$ there exists an $n \in \mathfrak{A}$ for which the equation \eqref{problem} has a solution in $k$. Then there exists an exponential polynomial  $a(n)$ with coefficients in a finite extension of $k$ that satisfies
\begin{equation} \label{solution}
a(n)^d+b_{d-1}(n)a(n)^{d-1}+ \dots +b_0(n)=0
\end{equation}
identically
\end{thm}
\begin{rem}
One could also try to prove something stronger than theorem \ref{main3}; namely that we have a solution to \eqref{problem} in the Hadamard ring as soon as we have solution for infinitely many $n$. A statement of this kind for the Hadamard quotient theorem is proved, with different methods, in \cite{finitezza} or in \cite{compositio} (the latter also deals in some cases with the root theorem). Some generalizations along the same lines are worked out in \cite{debrecen} and \cite{manuscripta}. Our present techniques do not allow us to obtain this stronger statement.
\end{rem}

The main theorem has a simple corollary, which deals with the case where the equation is not necessarily monic.
\begin{cor} \label{simple}
Let $k$ be a number field, $b_0, \dots, b_d \in \h{k}$, and suppose that for every $n$ the equation
\begin{equation} \label{not monic}
b_d(n)Y^d+b_{d-1}(n)Y^{d-1}+ \dots +b_0(n)=0
\end{equation} has a solution $a_n \in k$ for every $n$. Then there exists a finite extension $k'/k$ and two series $\sum a_1 (n) x^n, \sum a_2(n) x^n \in \h{k'}$ such that the sequence obtained as a component-wise quotient $a(n)=a_1(n)/a_2(n)$ (whenever defined) is a solution of \eqref{not monic}.
\end{cor}
To obtain the final form of our theorem we use for convenience a strengthening of the Hadamard quotient theorem, proved by Corvaja and Zannier in \cite{finitezza} (probably the result of \cite{VP88} suffices, but certainly leads to some difficulties). In that paper they use a form of the Subspace Theorem to prove that the conclusion of the Hadamard quotient theorem holds under the weaker hypothesis that the quotients $b(n)/c(n)$ lie in finitely generated ring for infinitely many $n$ (excluding some special cases). In section \ref{extension} we give a precise statement of a corollary of their theorem that we need. Combining this with corollary \ref{simple} we get our final result:
\begin{thm} \label{finale}
In the hypothesis of corollary \ref{simple} suppose moreover that the sequence of solutions $\{ a_n \}_{n \in \N}$ to \eqref{not monic} can be taken inside a finitely generated ring. Then there exists a finite extension $k'/k$ and a series $\sum a(n) x^n \in \h{k'}$ such that $a(n)$ is a solution of \eqref{not monic} for all $n$ such that $b_d(n) \neq 0$.
\end{thm}

\begin{rem}
A recent paper by Corvaja (\cite{Corvaja}) gives another perspective on these theorems. Corvaja restates our results in the context of actions of algebraic groups over algebraic varieties. The theory appears there because the entries of a power $A^n$ of a matrix are given by linear recurrences in $n$. In particular, Corvaja proves the following
\begin{thm*}[Corvaja]
Let $k$ be a number field and $G$ be a connected linear algebraic group, defined over $k$. Let $V$ be an affine algebraic variety and $\pi \colon V \mapsto G$ a finite map, both defined over $k$. Let $\Gamma \subset G(k)$ be a Zariski-dense semigroup. If $\Gamma$ is contained in the set $\pi(V(k))$, then there exists a connected component $V'$ of $V$ such that the restriction $\pi \res{V} \colon V' \mapsto G$ is an unramified cover. In particular $V'$ has the structure of an algebraic group over $k$.
\end{thm*}
As explained there, this can be seen as a geometric generalization of the Hilbert irreducibility theorem. Our result is used as a crucial starting point, giving the preceding assertion for the case where $\Gamma$ is cyclic.
\end{rem}

As we will see in the proofs, a central point of our argument is to guarantee that, given an absolutely irreducible polynomial $T(\ve{X},Y)$ over the number field $k$ (satisfying suitable conditons), we can find some suitable roots of unity $\{ \zeta_i \}$ such that the specialized polynomial $T(\zeta_1, \dots, \zeta_k, Y)$ remains irreducible over $k(\zeta_i)$. In the Master thesis \cite{Tesi} this was achieved with a reduction modulo some prime and an application of the Lang-Weil theorem. We give a description of this method in the appendix; although this approach is more complicated, it should be useful in other contexts. This step is simplified in the present proof by using a strong form of Hilbert irreducibility theorem for cyclotomic fields, obtained by Dvornicich and Zannier in \cite{StrongHilbert}; this work, in turn, is based on a result of Loxton (\cite{Loxton}), which bounds the number of addends necessary to write a cyclotomic integer $\alpha$ as a sum of roots of unity in terms of the maximum absolute value of the conjugates of $\alpha$ over $\Q$.

Before turning to the proofs we summarize here our notation.
\begin{description}[$a(n)$, $b(n)$, ,]
\item[$k$, $\ktl$] number fields
\item[$\mc{R}$] a ring of integers over a number field
\item[$\P$, $\mc{Q}$] prime ideals in  $\mc{R}$
\item[$\h{k}$] the Hadamard ring over the field $k$
\item[$k^c$] the maximal cyclotomic extension of a field $k$
\item[$a(n)$, $b(n)$] exponential polynomials, or the corresponding recurrence sequences
\item[$f$, $g$, $h$] polynomials
\item[$\ve{X}$] the vector of indeterminates $(X_1, \dots, X_r)$
\item[$\ve{a}$, $\ve{b}$] multiindices
\item[$\mathfrak{A}$, $\mathfrak{A}'$] arithmetic progressions
\item[$\gm$] the multiplicative group variety $GL_1$
\item[$\zeta$] some root of unity
\item[$\omega_n$] a primitive $n$-th root of unity
\end{description}
Note that we use a different symbol to distinguish between some generic root of unity and one of a fixed order.

\begin{ack}
We wish to thank Pietro Corvaja and Antonella Perucca for helpful comments.
\end{ack}

\section{Some reductions}

In the next sections we present the proof of theorem \ref{main}; in the present section we make some easy reductions, while the following section collects some techniques about the specialization of polynomials at roots of unity, which will be central in our argument.

The proof will be divided in several steps. The first two steps will fix some notation and make some reductions, while the crux of the arguments will appear from step \ref{crux} onwards. At the end of step \ref{polinomi S_D}, when we have fixed our notation, we present a brief sketch of how the proof will go on.

\begin{step}
Reduction to the case when the multiplicative subgroup generated by the $\beta_i$ inside $k^*$ is free.
\end{step}

We start with an easy lemma.
\begin{lemma} \label{riduzione}
In proving Theorem \ref{main3} it is possible to assume as well that the multiplicative subgroup $\Gamma < k^*$ generated by $\{ \beta_i \, | \, i=1, \dots, m\}$ is free.
\end{lemma}
\begin{proof}
Let $N$ be the order of the torsion part of $\Gamma$. Consider the exponential polynomials $b_{j,r}(n)=b_j(r+Nn)$, for some fixed $r$, $0 \leq r \leq N-1$; their roots are the $\beta_i^N$, so they generate a torsion-free group. Suppose that the theorem holds under the hypothesis of this lemma: we then get some exponential polynomials $a_r(n)$ such that
\begin{equation*}
a_r(n)^d+b_{d-1,r}(n)a_r(n)^{d-1}+ \dots +b_{0,r}(n)=0.
\end{equation*}
We may choose exponential polynomials $c_r(n)$ such that $c_r (Nn)=a_r(n)$. We remark that the exponential polynomial
\begin{equation*}
\theta(n)=\frac{1}{N} \sum_{i=1}^{N} \omega_N^n
\end{equation*}
takes the value $1$ for $N|n$ and $0$ otherwise. We define
\begin{equation*}
a(n)=\sum_{r=0}^{N-1} \theta(n-r) c_r(n-r).
\end{equation*}
In this way if $n=s+Nm$, with $0 \leq s \leq N-1$, we find $a(n)=a(s+Nm)= c_s(Nm)=a_s(m)$, and so equation \eqref{solution} is satisfied.
\end{proof}

We shall henceforth work under the additional hypothesis that $\Gamma$ is free. Having chosen a multiplicative basis $\gamma_1, \dots, \gamma_r$ we can write
\begin{equation*}
b_j(n)=f_j(n,\gamma_1^n, \dots, \gamma_r^n),
\end{equation*}
where the $f_j$ are rational function in $X_0, \dots, X_r$ of the special form
\begin{equation*}
f_j(X_0, \dots, X_r)=\frac{\widetilde{f_j}(X_0, \dots, X_r)}{X_1^{a_1} \cdots X_r^{a_r}},
\end{equation*}
$\widetilde{f}_j$ a polynomial. We call such a rational function a \emph{Laurent polynomial}; for all we need to do in this paper Laurent polynomials behave much like the classical ones. In particular the ring of Laurent polynomials is a localization of $\cl{k}[X_0, \dots, X_r]$, hence a $UFD$.

\begin{step}\label{polinomi S_D}
Reduction to the problem of proving that some equations have solution in a polynomial ring.
\end{step}

Consider the equation
\begin{equation} \label{pol prob}
Y^d+f_{d-1}(X_0, X_1^D \dots, X_r^D) Y^{d-1}+ \dots +f_{0}(X_0, X_1^D \dots, X_r^D)=0
\end{equation}
where we look for a solution $Y=Y(X_0,\dots, X_r)$ in the form of a Laurent polynomial in $X_0, \dots, X_r$. If \eqref{pol prob} has a solution the theorem is proved: it is sufficient to put
\begin{equation*}
a(n)=Y(n,\alpha_1^n, \dots, \alpha_r^n),
\end{equation*}
where $\alpha_i$ is a $D$-th root of $\gamma_i$. By construction \eqref{solution} holds.

\begin{rem}
As $n$ varies in $\N$, the $(r+1)$-uple $(n,\gamma_1^n, \dots, \gamma_r^n)$ describes a cyclic subsemigroup $C$ of $\mathbb{A}^1 \times \gm^r$. The equation \eqref{pol prob} defines a subvariety $V$ of $\mathbb{A}^1 \times \gm^r \times \mathbb{A}^1$; projection on the first $r+1$ coordinates gives a ramified covering of degree $d$
\begin{equation*}
\pi \colon V \to \mathbb{A}^1 \times \gm^r.
\end{equation*}

The hypothesis that equation \eqref{problem} has a solution for all $n$ can be rephrased saying that $C \subset \pi (V(k))$. The conclusion that we are trying to obtain is that for some $D$ there is a Laurent polynomial $Y(X_0,\dots, X_r)$ satisfying \eqref{pol prob}. Consider the unramified covering
\begin{equation*}
\begin{CD}
\rho_D \colon \mathbb{A}^1 \times \gm^r @>>> \mathbb{A}^1 \times \gm^r.
\\
(X_0,\dots, X_r) @>>> (X_0, X_1^D, \dots, X_r^D)
\end{CD}
\end{equation*}
This induces a cartesian diagram
\begin{equation*}
\begin{CD}
V' @>>> V
\\
@V{\pi'}VV @VV{\pi}V
\\
\mathbb{A}^1 \times \gm^r @>>{\rho_D}> \mathbb{A}^1 \times \gm^r,
\end{CD}
\end{equation*}
where $V'$ is the fibered product of $V$ and $\mathbb{A}^1 \times \gm^r$. The Laurent polynomial $Y(X_0,\dots, X_r)$ gives rise to a section $\tau \colon \mathbb{A}^1 \times \gm^r \to V'$ of $\pi'$; the existence of this section means that a component of $V'$ is a trivial covering of $\mathbb{A}^1 \times \gm^r$. This implies that some component of $V$ doesn't ramify over $\mathbb{A}^1 \times \gm^r$.

This is the point of view of \cite{StrongHilbert} (see theorem 1), of \cite{notePisa} (see the conjecure at p. $62$) and of \cite{Corvaja}, where this construction is generalized to ramified coverings of connected linear algebraic groups.
\end{rem}

To prove Theorem \ref{main2} we can thus assume that for each $D \geq 1$ the equation \eqref{pol prob} doesn't have a solution in the form of a Laurent polynomial. Gauss' lemma guarantees that the same equation doesn't have solutions in $\cl{k}(X_0,\dots,X_r)$. Define the Laurent polynomials
\begin{equation*}
S_D(X_0,\ve{X},Y)=Y^d+f_{d-1}(X_0, X_1^D \dots, X_r^D) Y^{d-1}+ \dots +f_{0}(X_0, X_1^D \dots, X_r^D);
\end{equation*}
our hypothesis is that these polynomials don't have linear factors in $Y$.

\begin{sketch}
The rest of the proof will be as follows. We consider $S_D$ for highly divisible values of $D$, we factorize it and work with one of the factors, call it $T$. We will be able to show that, since $\deg_Y T \geq 2$,
\begin{center}
\emph{there is some arithmetic progression $\mathfrak{A}$ such that for all $n \in \mathfrak{A}$ the specialization $T(n,\gamma_1^n, \dots, \gamma_r^n, Y)$ does not have roots in the base field.}
\end{center}
This is the main arithmetical point (it is almost the thesis of theorem \ref{main3}); it will be achieved in two steps.

First we show that the same property holds for most specializations of $T$ at roots of unity; namely if $(\zeta_0, \dots, \zeta_r)$ are generic roots of unity, then the specialized polynomial $T(\zeta_0, \dots, \zeta_r, Y)$ does not have roots in $k$. Actually we obtain the stronger result that it does not have solutions $\bmod \mc{Q}$ for some suitable ideal $\mc{Q}$ in the ring of integers of $k$. Hence in the next section we study a criterion for the irreducibility of the specialization of polynomials at roots of unity.

For the second step we use Chebotarev's theorem in order to choose roots of unity $\zeta_i$ that satisfy the congruences $\zeta_0 \equiv n$ and $\zeta_i \equiv \gamma_i^n$ $\pmod{\mc{Q}}$ whenever $n$ ranges in an arithmetic progression $\mathfrak{A}$. Combining these two steps we obtain the claim.

This takes already care of all the cases when $S_D$ is irreducible (and so equals $T$), for example the cyclotomic case treated in \cite{Z00}.

If $S_D$ is reducible, then  we make a change of variables, in order to restrict our exponential polynomials to the progression $\mathfrak{A}$. Then we repeat the same procedure with another factor of $S_D$, and so on. If in the process we end up with a linear factor, we are done; otherwise we end up with an arithmetic progression $\mathfrak{A}'$ such that \eqref{problem} does not have solution for $n \in \mathfrak{A}'$.
\end{sketch}

\section{Specialization of polynomials at roots of unity}

\begin{step} \label{crux}
A form of Hilbert irreducibility theorem for specializations at roots of unity.
\end{step}

We will now prove the following result about the specialization of Laurent polynomials, as a corollary of a work by Dvornicich and Zannier (\cite{StrongHilbert}):
\begin{prop} \label{Metodo Loxton}
Let $k$ be a number field and denote by $k^c$ its maximal cyclotomic extension. Let $f$ be a Laurent polynomial with coefficients in  $k^c$ and suppose that $f(X_1^{\ve{a}_1}, \dots, X_r^{\ve{a}_r}, Y)$ is irreducible over $k^c$ for every multiindex $\ve{a}$ with each $\ve{a}_i \leq \deg_Y f$. Then there exists a subvariety $W \subsetneq \gm^{r+1}$ such that if the $\zeta_i$ are roots of unity and $(\zeta_1, \dots, \zeta_r) \notin W$, the specialized polynomial $f(\zeta_1, \dots, \zeta_r,Y)$ is irreducible in $k^c[Y]$.
\end{prop}
We shall make use of the following result from \cite[\S 1.2, Lemma 2]{Schinzel}
\begin{prop} \label{schinzel}
Let $K$ be a field and $f \in K[\ve{X},Y]$; there exist a polynomial $g \in K[\ve{X},Y]$ and a non-zero polynomial $g_1 \in K[\ve{X}]$ with the following property. Suppose $x_1, \dots, x_r$ lie in some extension $L$ of $K$ and $g_1(x_1, \dots, x_r) \neq 0$; then $f(x_1, \dots, x_r, Y)$ is reducible in $L[Y]$ if, and only if, $g(x_1, \dots, x_r, Y)$ has a root in $L$.
\end{prop}
Actually the proposition is stated there for polynomials, but it is easy to derive the conclusion for Laurent polynomials as well. To prove proposition \ref{Metodo Loxton} we will also need the following
\begin{prop} \label{fattori}
Let $A(X_1, \dots, X_r, Y)$ be a Laurent polynomial with coefficients in some field $k$, and suppose that $A(\ve{X}^{\ve{a}_1}, \dots, \ve{X}^{\ve{a}_r}, Y)$ is reducible over $k$ for some multiindices $\ve{a}_1, \dots , \ve{a}_r \in \Z^r$. Suppose moreover that the $\ve{a}_i$ are linearly independent. Then there is a $m \leq \deg_Y A$ such that $A(X_1^m, \dots, X_r^m, Y)$ is reducible.
\end{prop}
\begin{proof}
It is easy to see that any lattice $\mathcal{L}$ inside $\Z^r$ contains a sublattice of the form
\begin{equation*}
\left\langle (M,0,\dots,0), (0,M,0,\dots,0), \dots, (0,\dots,0,M) \right\rangle,
\end{equation*}
where $M$ is the discriminant of $\mathcal{L}$. In fact if $B(\mathcal{L})$ is a matrix whose columns fom a basis of $\mathcal{L}$ and $B' (\mathcal{L})$ is the cofactors matrix, then $B'(\mathcal{L}) \cdot B(\mathcal{L}) = M \id$.

Moreover if the $\ve{b}_j$ form a sublattice of the lattice spanned by the $\ve{a}_j$, then by substitution we obtain that $A(\ve{X}^{\ve{b}_1}, \dots, \ve{X}^{\ve{b}_r}, Y)$ is reducible too. Combining these facts we can assume that we have a factorization
\begin{equation*}
A(X_1^M, \dots, X_r^M, Y) = A_1(X_1, \dots, X_r, Y) \cdots A_m(X_1, \dots, X_r, Y)
\end{equation*}
for some $M \in \N$. We get an action of $\left( \Z/M\Z \right)^r$ on the set $\{ A_1, \dots, A_m \}$ of factors by letting
\begin{equation*}
(a_1, \dots, a_r).A_i (X_1, \dots, X_r, Y) = A_i (\omega_M^{a_1} X_1, \dots, \omega_M^{a_r} X_r, Y).
\end{equation*}
The index of the stabilizer of the factor $A_1$ is $m'=\# \Orb (A_1) \leq m$, hence this stabilizer contains a subgroup of the form
\begin{equation*}
k_1\Z/M\Z \times \dots \times k_r\Z/M\Z,
\end{equation*}
where $k_i$ divides $m'$. This means that each monomial in $A_1$ involves the variable $X_i$ at a power multiple of $M/k_i$, which in turn is multiple of $M/m'$; hence we can write $A_1(X_1, \dots, X_r, Y) = A'_1(X_1^{M/m'}, \dots, X_r^{M/m'}, Y)$. The same holds for the complementary factor. But this implies that $A(X_1^{m'}, \dots, X_r^{m'}, Y)$ is already reducible, and by construction $m' \leq \deg_Y A$
\end{proof}

\begin{proof}[Proof of proposition \ref{Metodo Loxton}]
By contradiction. Assume that there exists a set $Z$ of roots of unity, Zariski dense in $\gm$, such that $f(\zeta_1, \dots, \zeta_r,Y)$ is reducible for each choice of $(\zeta_1, \dots, \zeta_r) \in Z$. With the notation of proposition \ref{schinzel}, it is not restrictive to suppose that for $(\zeta_1, \dots, \zeta_r) \in Z$ we have $g_1(\zeta_1, \dots, \zeta_r) \neq 0$; then proposition \ref{schinzel} guarantees that $g(\zeta_1, \dots, \zeta_r, Y)$ has a root in $k^c$. If $g$ is reducible, there is at least one of his irreducible factors $g_2$ such that the subset of $Z$ for which $g_2(\zeta_1, \dots, \zeta_r, Y)$ has a root in $k^c$ is still dense; we replace $Z$ by this smaller subset.

We apply theorem $1$ of \cite{StrongHilbert} with $V$ the zero locus of $g_2$ inside $\gm^{r+1}$ and $\pi \colon V \mapsto \gm^r$ the projection on the $X$ coordinates. The hypothesis of the theorem require that the subset $J$ of $V$ consisting of those elements mapping to roots of unity is dense in $V$. By construction we know that $\pi (J) \supset Z$, so it is dense in $\gm^r$. It follows that $\dim \overline{J} \geq r$, so $J$ is actually dense in $V$ by irreducibility.

The theorem gives us a lot of information. First, the closure of $\pi(V)$ has the form $\zeta T$, where $T$ is a subtorus of $\gm^r$, and $\zeta$ is torsion. In our case $T$ equals $\gm^r$, since we already know that $\pi(V)$ is dense. Moreover we get the existence of an isogeny $\mu \colon T \mapsto T$ and a rational map $\rho \colon T \dashrightarrow V$,  defined over $k^c$, such that $\pi \circ \rho = \zeta \cdot \mu$.

In our situation we can assume that $\zeta=1$, since $T$ is the whole $\gm^r$. Moreover it is well known that the isogeny $\mu \colon \gm^r \mapsto \gm^r$ must be of the form
\begin{equation*}
(x_1, \dots, x_r) \mapsto (\ve{x}^{\ve{a}_1}, \dots, \ve{x}^{\ve{a}_r})
\end{equation*}
for suitable linearly independent multiindices $\ve{a}_i$. The rational map $\rho$ can be written as $\big( R_1(X_1, \dots, X_r), \dots, R_{r+1}(X_1, \dots, X_r) \big)$, where the $R_i \in k^c(X_1, \dots, X_r)$. The fact that $\rho$ takes values in $V$ can be translated saying that
\begin{equation*}
g_2 \big(R_1(X_1, \dots, X_r), \dots, R_{r+1}(X_1, \dots, X_r) \big) =0.
\end{equation*}
The fact that it is, up to isogeny, a section of $\pi$ means that $R_i(X_1, \dots, X_r)=\ve{X}^{\ve{a}_i}$ for $i = 1, \dots, r$. So \emph{a fortiori}
\begin{equation*}
g \big(\ve{X}^{\ve{a}_1}, \dots, \ve{X}^{\ve{a}_r}, R_{r+1}(X_1, \dots, X_r) \big) =0.
\end{equation*}
This means that $g (\ve{X}^{\ve{a}_1}, \dots, \ve{X}^{\ve{a}_r},  Y)$ has a root in $k^c(X_1, \dots, X_r)$; by proposition \ref{schinzel} again we obtain that $f(\ve{X}^{\ve{a}_1}, \dots, \ve{X}^{\ve{a}_r},  Y)$ is reducible over $k^c$. Proposition \ref{fattori} now allows us to conclude.
\end{proof}

\begin{step}
The irreducibility properties of our polynomials.
\end{step}

We don't know very much about the irreducibility of the Laurent polynomials $S_D$, but let us vary $D$, making it more and more divisible. The number of factors will stabilize to a number less than $\deg_Y g$, since $g$ is monic in the $Y$ variable. So there is a $D_0$ such that if $S_{D_0}$ factors as
\begin{equation*}
S_{D_0}(X_0, \ve{X},Y) = T_1(X_0, \ve{X},Y) \cdots T_l(X_0, \ve{X},Y),
\end{equation*}
then for every $M \in \N$
\begin{equation*}
S_{M D_0}(X_0, \ve{X},Y) = T_1(X_0, X_1^M, \dots, X_r^M, Y) \cdots T_l(X_0, X_1^M, \dots, X_r^M, Y)
\end{equation*}
will also be a decomposition into prime factors. Our hypothesis in step \ref{polinomi S_D} amounts to saying that $\deg_Y T_i \geq 2$ for each $i=1,\dots l$.

It is not restrictive to assume that $D_0=1$, as we shall do from now on. In fact multiplying by $D_0$ the terms of an arithmetic progression yields another arithmetic progression (see also step \ref{conclusione}). We now want to specialize the first variable $X_0$ in such a way to preserve irreducibility. By Hilbert irreducibility theorem (\cite[\S 4.4]{Schinzel}) we can find some $\theta \in k$ such that each factor $T_j(\theta, X_1^m, \dots, X_r^m, Y)$ remains irreducible for $m \leq \deg_Y T_j$. Proposition \ref{fattori} guarantees that $T_j(\theta, \ve{X}^{\ve{a_r}}, \dots, \ve{X}^{\ve{a_r}}, Y)$ will be irreducible for each choice of linearly independent multiindices $\ve{a}_j$.

\section{Proof of the main theorem}

\begin{step}
Some irreducible factor $T$ of $S_D$ admits an irreducible specialization at roots of unity.
\end{step}

We choose a rational prime $\beta$ multiplicatively independent from the $\gamma_i$s, and put $\delta_i=\gamma_i \beta^{k_i}$, for some integers $k_i$ which we shall choose later. The following lemma is proved by Zannier in \cite{Z00} (this is where we make use of the fact that the multiplicative group $\Gamma$ is free).
\begin{lemma}[\cite{Z00}] \label{lemmazannier}
There exists a number $L$ such that, whenever we take $M \geq 1$ and a prime $\ell>L$, then $\beta^M$ doesn't belong to the multiplicative group generated by the $\delta_i$ and by $((k^c)^*)^{\ell M}$. The number $L$ depends on $k$, $\beta$ and $\gamma_i$, but it doesn't depend on the $k_i$.
\end{lemma}

We fix once and for all a natural number $L$ greater than $\deg_Y g$ and big enough for the preceding Lemma to hold. Consequently we choose $D$ divisible by each prime factor less than $L$ and big enough, so that the inclusion $\Q^c \cap k \subset \Q(\omega_D)$ holds. The latter choice will guarantee that for each $s \geq 1$, $\Q(\omega_{sD})/\Q(\omega_{D})$ and $k (\omega_{D}) / \Q (\omega_{D})$ are linearly disjoint extensions.

Now we fix some factor $T$, say $T_1$, of $S$; we will work with this polynomial until the last step. Let us put
\begin{equation*}
\widetilde{T} (X_1, \dots, X_r, Y)=T(\theta, X_1^{D}, \dots, X_r^{D}, Y),
\end{equation*}
where $\theta$ is defined at the end of the previous section.

\begin{lemma}
Let $W \subsetneq \gm^r$ be an algebraic subvariety of a torus, defined over $k$, and fix a natural number $M$. Then there exist roots of unity $\zeta_1, \dots, \zeta_r$ such that:
\begin{enumerate}[i)]
\item
$(\zeta_1, \dots, \zeta_r) \notin W(k^c)$
\item
the order of each $\zeta_j$ is not multiple of a prime less than $M$.
\end{enumerate}
\end{lemma}

\begin{proof}
Let $S$ be the set of roots of unity whose order is not multiple of a prime less than $M$. Since $S$ is infinite, it is dense in $\gm$, and this is the thesis in the case $r=1$. In the general case $S^r$ is dense in $\gm^r$.
\end{proof}


The preceding lemma, together with proposition \ref{Metodo Loxton}, allows us to fix roots of unity $\zeta_1, \dots, \zeta_r$  such that the multiplicative order of $\zeta_j$ is not multiple of a prime smaller than $L$, and at the same time
\begin{equation} \label{acca}
h(Y)=\widetilde{T} (\zeta_1, \dots, \zeta_r, Y)
\end{equation}
remains irreducible over $k$.

\begin{step}
The specialized polynomial $h$ has no roots, even modulo some suitable primes.
\end{step}

We start with a lemma.
\begin{lemma} \label{scelta di p}
There exist infinitely many primes of the form $p=1+Dm$ such that
\begin{enumerate}[i)]
\item
every prime factor of $m$ is greater than $L$
\item
we can write $\zeta_i = \omega_{p-1}^{k_i}$ for suitable integers $k_i$.
\end{enumerate}
\end{lemma}

\begin{proof}
This is an easy consequence of Dirichlet's theorem on the existence of primes in arithmetic progressions. Let $s$ be the lowest common multiple of the orders of $\zeta_0, \dots, \zeta_r$. We need a prime $p$ satisfying the following congruences:
\begin{equation} \label{congruenze Dirichlet}
\begin{cases}
p \equiv 1 \pmod{D s}
\\
p \nequiv 1 \pmod{D \ell} \text{ for each prime } \ell \leq L.
\end{cases}
\end{equation}
Indeed the first congruence guarantees that $p$ can be written in the form $1+Dm$ for some $m$, and that $p-1$ is multiple of the order of every root of unity $\zeta_i$, while the second condition implies that $m$ doesn't have any prime factor smaller than $L$. Thanks to the chinese remainder theorem and Dirichlet's theorem we find a prime solution to \eqref{congruenze Dirichlet}.
\end{proof}


The preceding lemma allows to fix the numbers $k_i$ mentioned at the beginning of the section. We define
\begin{equation}
\begin{split}
\ktl &= k(\omega_{p-1})
\\
E &= k( \omega_{p-1}, \beta^{1/m}, \delta_1^{1/m}, \dots, \delta_r^{1/m}).
\end{split}
\end{equation}

\begin{lemma}
The polynomial $h$ defined in \eqref{acca} remains irreducible in $E$.
\end{lemma}
\begin{proof}
Assume this is not the case, and factor $h$ as $h=h_{1} h_{2}$ where $0<d_i=\deg h_{i}<\deg h$. Let $E'$ be obtained by adding a root of $h_{1}$ to $E$. By Kummer theory we know that $[E:\ktl]$ divides a power of $m$, so $[E':\ktl]$ divides $d_1$ times a power of $m$. On the other hand, by construction $h$ admits a root in $E'$, so $\deg h$ divides $[E':\ktl]$; this is impossible since each prime factor of $m$ is $>L \geq d$.
\end{proof}

\begin{lemma}
\begin{equation}
\label{grado beta} [E:k(\omega_{p-1}, \delta_1^{1/m}, \dots, \delta_r^{1/m})]=m.
\end{equation}
\end{lemma}
\begin{proof}
If the degree were lower, it would be a proper divisor of $m$, again by Kummer theory. Take a prime $\ell$ such that this degree divides $m/\ell$. We can apply Kummer theory to the field $\ktl$: the two groups
\begin{align*}
\Delta= & \langle (\ktl^*)^m, \beta, \delta_1,\dots, \delta_r \rangle
\\
\Delta^{\prime}= & \langle(\ktl^*)^m, \beta^{\ell}, \delta_1,\dots, \delta_r \rangle
\end{align*}
define the same extension $E/\ktl$, so they coincide. In particular we can express
\begin{equation*}
\beta =\alpha^m \beta^{\ell a_0} \delta_1^{a_1} \cdots \delta_r^{a_r}
\end{equation*}
for some $\alpha \in \ktl$. But this contradicts lemma \ref{lemmazannier} with $M=1$ (note that $\ell>L$).
\end{proof}

Since the extension $E/k(\omega_{p-1}, \delta_1^{1/m}, \dots, \delta_r^{1/m})$ is cyclic we can take a generator $\tau$ of its Galois group.
\begin{lemma} \label{galois sposta}
Call $E'$ the splitting field of $h(Y)$ over $E$. There exists $\xi \in G=\Gal (E',\ktl)$ such that:
\begin{enumerate}[i)]
\item $\xi\res{E}=\tau$
\item if $y$ is a root of $h$, then $\xi (y) \neq y$.
\end{enumerate}
\end{lemma}

\begin{proof}
We first show the existence of some $\sigma \in G$ satisfying $ii)$. This amounts to prove that the union of the stabilizers of the roots of $h$ is not all of $G$. By the irreducibility of $h$ these stabilizers are conjugate subgroups. Let $H$ be one of them; then there are at most $|G|/|H|$ stabilizers, each one of order $|H|$, so the union can't be all of $G$ (they all contain the identity).

Let $\widetilde{\sigma}=\sigma\res{E} \in \Gal (E,\ktl)$, and $\varphi=\widetilde{\sigma}^{-1} \tau$. We only need to extend $\varphi$ to $E'$ in such a way that $\varphi (y)=y$ for every root $y$ of $h$. If we call $F$ the splitting field of $h$ over $\ktl$, so that $E'=E F$, we reduce to the problem of verifying that $E$ and $F$ are linearly disjoint over $\ktl$. This follows by comparison of the degrees: $[F:\ktl]$ divides $d!$, while $[E:\ktl]$ divides some power of $m$, and each prime factor of $m$ is $>L \geq d$.
\end{proof}

Let $\mc{R}$ be the ring of integers of $\ktl$. By Chebotarev theorem we get a positive density set of primes $\mc{Q}$ of $\mc{R}$ whose Frobenius verifies $\phi(\mc{Q}'|\mc{Q})=\xi$ in $E'$, for some prime $\mc{Q}'$ over $\mc{Q}$. We don't affect the density if we ask that $\mc{Q}$ has no inertia over the rationals. Subject to these conditions, we take a big prime $\mc{Q}$ at which the reductions of $\beta$, $\gamma_j$ and $f$ are defined, and call $\mc{Q'}$ a prime over it such that $\phi(\mc{Q'}|\mc{Q})=\xi$.

\begin{lemma} \label{noradici}
If $\mc{Q}$ has big enough norm, then the congruence $h(Y)\equiv 0 \pmod{\mc{Q}}$ has no solutions.
\end{lemma}
\begin{proof}
First we remark that if $\mc{Q}$ is big enough, $h$ has distinct roots $\bmod{\, \mc{Q}}$. Let $y$ be one of such roots: by lemma \ref{galois sposta} we know that $\xi(y)\neq y$. If $\mc{Q}$ has big enough norm and $\mc{Q}'$ is above $\mc{Q}$, then $\mc{Q}'$ is not a prime factor of the number $\xi(y)-y$ for any root $y$ of $h$. Hence for every root $y$ we have $\xi(y)\nequiv y \pmod{\mc{Q'}}$.

This means that the Frobenius of $\F_{q}=\mc{R}/\mc{Q}$ doesn't fix the class $\cl{y} \in \cl{\F_{q}}$, that is, $h$ has no roots $\bmod{\, \mc{Q}}$.
\end{proof}

\begin{step}
If $n$ is chosen in a suitable arithmetic progression, then the polynomial $T(n,\gamma_1^n, \dots, \gamma_r^n, Y)$ has no roots in the base field.
\end{step}

\begin{lemma}
There exists an arithmetic progression $\mathfrak{A}$ such that if $n \in \mathfrak{A}$, then for each $j=1, \dots, r$
\begin{equation*}
\begin{cases}
n \equiv \theta &\pmod{\mc{Q}}
\\
\gamma_j^n \equiv \omega_{m}^{k_j} &\pmod{\mc{Q}}.
\end{cases}
\end{equation*}
\end{lemma}
\begin{proof}
By our choices we  know that $q=N_{\Q}^{\ktl}(\mc{Q})$ splits completely in $\Q(\omega_{p-1})$, so we deduce that $q \equiv 1 \pmod{p-1}$, and in particular $m|q-1$. Moreover $\xi$ fixes each $\delta_j^{1/m}$, so $\delta_j$ is a $m$-th power $\bmod{\, \mc{Q}}$, hence
\begin{equation*}
\delta_j^{\frac{q-1}{m}}\equiv 1 \pmod{\mc{Q}}.
\end{equation*}
Similarly $\xi(\beta)=\omega_m^a \beta$ for some $a$, which is coprime with $m$ by (\ref{grado beta}), so
\begin{equation*}
\beta^{\frac{q-1}{m}}\equiv \omega_m^a \pmod{\mc{Q}}.
\end{equation*}
Putting the two relations together we deduce
\begin{equation*}
\gamma_j^{\frac{q-1}{m}}\equiv \omega_m^{a k_j} \pmod{\mc{Q}}.
\end{equation*}
Calling $b$ the inverse of $a$ $\bmod{\, m}$ we find that
\begin{equation*}
\gamma_j^{b\frac{q-1}{m}}\equiv \omega_m^{k_j} \pmod{\mc{Q}}.
\end{equation*}
Moreover we can take $c \in \N$ satisfying $c \equiv \theta \pmod{\mc{Q}}$. If $n \in \N$ is a solution of the congruences
\begin{equation*}
n \equiv c \pmod{q}, \qquad n\equiv b\frac{q-1}{m} \pmod{q-1},
\end{equation*}
then we have the relations
\begin{equation} \label{congruenze}
\begin{cases}
n \equiv \theta &\pmod{\mc{Q}}
\\
\gamma_j^n \equiv \omega_{m}^{k_j} &\pmod{\mc{Q}} \quad j=1,\dots,r.
\end{cases}
\end{equation}

\end{proof}

\begin{lemma}
Assume that $n$ is taken in the arithmetic progression $\mathfrak{A}$. Then the polynomial $T(n,\gamma_1^n, \dots, \gamma_r^n, Y)$ has no roots in $k$.
\end{lemma}
\begin{proof}
The conditions \eqref{congruenze} imply that
\begin{equation*}
T(n,\gamma_1^n,\dots,\gamma_r^n,Y) \equiv T(\theta, \omega_{p-1}^{k_1}, \dots, \omega_{p-1}^{k_r}, Y) \equiv h(Y) \pmod{\mc{Q}}.
\end{equation*}
If $T(n,\gamma_1^n,\dots,\gamma_r^n,Y)$ had a root in $k$, then $h$ would have a root $\bmod{\, \mc{Q}}$, which is excluded by lemma \ref{noradici}
\end{proof}

\begin{step}\label{conclusione}
Conclusion of the proof of theorem \ref{main3}.
\end{step}

Now if $T_1$ is the only factor of $S$, we are done. Otherwise we proceed in the following way. First, we describe the arithmetic progression
\begin{equation*}
\mathfrak{A}=\{ a n +b, \,  n \in \N \}
\end{equation*}
for suitable $a,b \in \N$. Next, we operate the substitution
\begin{equation*}
T_i'(X_0, X_1, \dots, X_r, Y)=T_i(a X_0+b,\gamma_1^b X_1^a, \dots,\gamma_r^b X_r^a, Y).
\end{equation*}
The $T_i'$ may not be irreducible anymore, but after further factorization and relabeling we assume that $T_2'$ is irreducible. If $T_2'$ has degree greater than one, we call it $T$ and repeat the whole procedure on and on. Eventually one of the following cases will happen:
\begin{enumerate}[i)]
\item
We get an arithmetic progression $\mathfrak{A}'=\{ a' n +b', \,  n \in \N \}$ and a degree one (in the last variable) factor of $S(a' X_0+b',\gamma_1^{b'} X_1^{a'}, \dots,\gamma_r^{b'} X_r^{a'}, Y)$, say $Y-Y(X_0, \dots, X_n)$. In this case let us take $\alpha_i$ such that $\alpha_i^{a'}=\gamma_i$. Put $a(n)=Y(n/a', \alpha_1^{n} \dots, \alpha_r^{n})$; then $a(n)$ is an exponential polynomial, and the relation
\begin{equation*}
S(a' X_0+b',\gamma_1^{b'} X_1^{a'}, \dots,\gamma_r^{b'} X_r^{a'}, Y(X_0, \dots, X_n))=0
\end{equation*}
gives, for $X_0=n/a'$ and $X_i=\alpha_i^n$,
\begin{equation*}
S(n+b',\gamma_1^{n+b'}, \dots,\gamma_r^{n+b'}, a(n))=0,
\end{equation*}
that is
\begin{equation*}
a(n)^d+b_{d-1}(n+b)a(n)^{d-1}+ \dots +b_0(n+b)=0,
\end{equation*}
so we have a solution of the original equation in the Hadamard ring.
\item
We never get a degree one factor. In this case, after at most $d/2$ steps we end with an arithmetic progression $\mathfrak{A}'=\{ a' n +b', \,  n \in \N \}$ such that \eqref{problem} has no solution in $k$ for $n \in \mathfrak{A}'$, which is the thesis.
\end{enumerate}
\qed

\section{Proof of the remaining assertions} \label{extension}

The aim of the present section is to prove corollary \ref{simple} and theorem \ref{finale}, which deal with not necessarily monic equations.

\begin{proof}[Proof of corollary \ref{simple}]
Multiplying \eqref{not monic} by $b_d(n)^{d-1}$ and putting $Z=b_d(n) Y$ we obtain the equation
\begin{equation*}
Z^d+b_{d-1}(n)Z^{d-1} \dots + b_0(n)b_d(n)^{d-1}=0;
\end{equation*}
this has a solution $a_2(n) \in \h{k'}$ for some finite extension $k'/k$, thanks to theorem \ref{main}. Putting $a_1(n)=b_d(n)$ we get the thesis.
\end{proof}

Preliminary to the proof of theorem \ref{finale} we cite a stronger form of the Hadamard quotient theorem, due to Corvaja and Zannier (\cite[cor. 2]{finitezza})
\begin{thm*}
Let $k$ be a number field and $\mc{R} \subset k$ a finitely generated ring. Let $\sum b(n) x^n,  \sum c(n) x^n \in \h{k}$ and assume that their roots generate a torsion-free group. Then either $b(n)/c(n)$ is a recurrence sequence or the set of natural numbers $n$ for which $b(n)/c(n) \in \mc{R}$ has zero density.
\end{thm*}

We will also need the Skolem-Mahler-Lech theorem (see \cite{BetterKnown}).
\begin{thm*}[Skolem, Mahler, Lech]
Let $K$ be a field of characteristic $0$ and let $a(n)$ be a linear recurrence over $K$. Then the zero set of $a$
\begin{equation*}
\{ n \in \N \, | \, a(n)=0 \}
\end{equation*}
is the union of a finite set with a finite number of complete arithmetic progression.
\end{thm*}
By a \emph{complete} arithmetic progression we mean a set of the form $\{ ak+b \, | \, k \in \N \}$ for some $a \in \N$, $b \in \{0, \dots, a-1 \}$; for example $\{5,8,11,\dots \}$ is not complete ($2$ is missing).

\begin{proof}[Proof of theorem \ref{finale}]
By corollary \ref{simple} we know that we can find two recurrence sequences $\{ a_1(n) \}$ and $\{ a_2(n) \}$ such that $a_2(n)/a_1(n)$ satisfies equation \eqref{not monic} for every $n$ such that the quotient is defined. We can argue as in lemma \ref{riduzione} to restrict ourselves to the case where the roots of $a_1(n)$ and $a_2(n)$ generate a torsion-free group, call it $G$. Let us call $A$ the ring of the recurrence sequences with roots in $G$. $A$ is isomorphic to a polynomial ring over $k$, in particular it is a unique factorization domain. We can divide both $a_1$ and $a_2$ by their greatest common divisor in $A$, so we shall assume that $a_1$ and $a_2$ are relatively prime.

At first suppose that $b_d$ never vanishes; then the same holds true for $a_1$, which divides $b_d$. In particular the quotient $a_2(n)/a_1(n)$ is always defined. The polynomial
\begin{equation*}
b_d(n)Y^d+b_{d-1}(n)Y^{d-1}+ \dots +b_0(n)
\end{equation*}
is divisible by $a_1(n)Y-a_2(n)$ in $K[Y]$, where $K$ is the field of fractions of $A$; by Gauss' lemma the same is true in $A[Y]$. So we have a factorization of the original equation as
\begin{equation*}
\left( a_1(n) Y - a_2(n) \right) \left( c_{d-1}(n) Y^{d-1}+c_{d-2}(n) Y^{d-2}+ \dots c_0 (n) \right) =0,
\end{equation*}
for suitable recurrence sequences $c_i(n)$. By induction on the degree, we know that either the equation
\begin{equation} \label{induzione}
c_{d-1}(n) Y^{d-1}+c_{d-2}(n) Y^{d-2}+ \dots c_0 (n)=0
\end{equation}
has a solution in some Hadamard ring (in which case we are done), or it is not solvable in the field for $n$ in some arithmetic progression $\mathfrak{A}$. But then we must have $\widetilde{a}(n)=a_2(n)/a_1(n)$ for $n \in \mathfrak{A}$; by the theorem of Corvaja and Zannier the quotient of $a_2(n)$ by $a_1(n)$ is then a recurrence sequence itself.

Now consider the general case. By the theorem of Skolem-Mahler-Lech we know that set zero set of $b_d$ is a union of a finite number of elements and a finite number of complete arithmetic progressions. Since we are working in $\h{k}$ we can disregard the finite number of terms; so we can assume that there is an $m \in \N$ and some numbers $n_1, \dots, n_r \in \{ 0, \dots, m-1 \}$ such that $b_d(n) =0$ if, and only if, $n \equiv n_i \pmod{m}$ for some $i$.

Fix a number $c \in \{ 0, \dots, m-1 \}$ different from all the $n_i$, and consider the equation
\begin{equation*}
b_d(c+ nm)Y^d+b_{d-1}(c +nm)Y^{d-1}+ \dots +b_0(c+nm)=0.
\end{equation*}
The coefficients $b_i(c+ nm)$ are linear recurrences in $n$ (up to a finite number of terms), and by construction $b_d(c+ nm)$ never vanishes. By the first part of the proof we can find a series $\sum a_c(n) x^n \in \h{k'}$ such that $a_c(n)$ satisfies the equation for all $n$. For $c=n_i$ we can choose any linear recurrence $a_c$, for example put $a_c(n)=0$ for all $n$.

As we have seen in the proof of lemma \ref{riduzione}, the exponential polynomial
\begin{equation*}
\theta(n)=\frac{1}{m} \sum_{i=1}^{m} \omega_m^n
\end{equation*}
takes the value $1$ for $m|n$ and $0$ otherwise. Choose exponential polynomials $a_c'(n)$ such that $a_c'(mn)=a_c(n)$. We define
\begin{equation*}
a(n)=\sum_{r=0}^{m-1} \theta(n-r) a_r'(n-r).
\end{equation*}
By construction $a(c+ nm)=a_c(n)$  for all $c =0, \dots, m-1$, so $a(n)$ satisfies equation \eqref{not monic} whenever $b_d(n) \neq 0$.
\end{proof}

\section{A different approach to the proof} \label{appendice}

In this appendix we discuss a different approach to the proof, as outlined in \cite{Tesi}. The method described here is more similar to the original article \cite{Z00}, but some new difficulties arise with respect to the case of cyclotomic equations, since in the general case we don't have Kummer theory available. As we have seen, one of the main points in the proof proposition \ref{Metodo Loxton}: namely we have to guarantee that the polynomial $h(Y)$, obtained by specialization of a factor $T$ of $S$ at roots of unity, remains irreducible, knowing that we can assume $T$ absolutely irreducible.

In what follows we a describe a different way to prove this. The notation is the same as in the proof of the main theorem \ref{main}. Since this approach is not fully developed, some detail is missing. We believe anyway that this method may prove itself useful to solve similar problems, where the alternative way doesn't work.

Given the absolutely irreducible polynomial $T$, one can construct another absolutely irreducible polynomial $\widetilde{T}$ in the following way. We look at $T$ as a polynomial in the $Y$ variable over $k(X_0,\dots,X_r)$, take some root $Y_0$, and denote by $L$ the normal closure of $k(X_0,\dots,X_r,Y_0)$ over $k(X_0,\dots,X_r)$. Since $char (k) =0$ we can write $L=k(X_0,\dots,X_r,Y_1)$ for a suitable $Y_1 \in L$; we set $\widetilde{T}$ to be the minimal polynomial of $Y_1$ over $k(X_0,\dots,X_r)$. By construction, whenever a specialization of $\widetilde{T}$ has a root inside some field, the specialization of $T$ at the same values has $\deg_Y T$ roots (maybe repeated) in the same field.

Then we make use of the estimates given by the Lang-Weil theorem (\cite{Lang-Weil}) to obtain the following proposition. This method goes back to Eichler, S. D. Cohen and others (see for example \cite{field}); we give a proof of the result that we use, since later we will want to point out some possible modifications.
\begin{prop} \label{lw}
Let $k$ be a number field with $\R$ as its ring of integers, and suppose that $T, \widetilde{T} \in \R[X_0, \dots, X_r,Y]$ are as above. Then for every prime $\P$ of $\R$ of big enough norm we can find some $(r+1)$-uple $(x_0, \dots, x_r) \in \F_q := \R/\P$ such that the equation $f(x_0, \dots, x_r, Y) \equiv 0$ has no solution in $\F_q$. Moreover we can assume that no $x_i$ is $0$ in $\F_q$.
\end{prop}

\begin{proof}
By a theorem of Owstrowski we know that the reduction of $T$ and $\widetilde{T}$ modulo $\P$ remain absolutely irreducible for $|\P|$ large. Let $\F_q$ be the residue field at $\P$ and call $N(q)$ the number of solutions to $\widetilde{T} \equiv 0$ in $\F_q^{r+2}$. Applying Lang-Weil we deduce that
\begin{equation*}
N(q)=q^{r+1}+O(q^{r+1/2}).
\end{equation*}
We know that if $(x_0,\dots , x_r) \in \F_q^{r+1}$ is such that $\widetilde{T}(x_0,\dots ,x_r, Y)$ has at least a solution, then $f(x_0,\dots ,x_r, Y)$ will have exactly $d_T = \deg_Y T$ solutions. Actually we should take care of repeated roots, but those will account only for a term $O(q^{r+1/2})$ in our estimates. The number of such $(r+1)$-uples is at least
\begin{equation*}
\frac{N(q)}{d_{\widetilde{T}}} \geq \frac{q^{r+1}}{d_{\widetilde{T}}}+O(q^{r+1/2}),
\end{equation*}
so we get at least $\frac{d_T}{d_{\widetilde{T}}}q^{r+1}+O(q^{r+1/2})$ solutions for $T$. Let us call $M(q)$ the number of the solutions for $T$ that we haven't counted yet. We can apply Lang-Weil, this time to $T$, and get
\begin{equation*}
\frac{d_{\widetilde{T}}}{d_T} q^{r+1}+M(q)=q^{r+1}+O(q^{r+1/2}),
\end{equation*}
which gives
\begin{equation*}
M(q) \leq \left( 1-\frac{d_{\widetilde{T}}}{d_T} \right) q^{r+1}+O(q^{r+1/2}).
\end{equation*}
It follows that the number of $(r+1)$-uples $(x_0,\dots, x_r)$ for which $T$ has at least a solution can be estimated by
\begin{equation} \label{stima}
\frac{q^{r+1}}{d_{\widetilde{T}}}+M(q)+O(q^{r+1/2}) \leq q^{r+1} \left( 1-\frac{d_T -1}{d_{\widetilde{T}}} \right)
+O(q^{r+1/2}).
\end{equation}
This is less than $q^{r+1}$ when $q$ is big enough, so the conclusion follows. To get the sharper statement it is enough to observe that the number of $(r+1)$-uples $(x_0,\dots, x_r)$ for which at least one of the $x_i$ is $0$ is trivially $O(q^r)$.
\end{proof}

Choose a prime $\P$ satisfying the conclusion of the preceding lemma, and without inertia over the rationals, so that $\R/\P = \F_p$, $p$ a prime. In the field $\ktl=k(\omega_{p-1})$ we can take representatives $(\omega_{p-1}^{a_0}, \dots, \omega_{p-1}^{a_r})$ for $(x_0, \dots, x_r)$ which are roots of unity, so we can conclude that the specialized polynomial $T(\omega_{p-1}^{a_0}, \dots, \omega_{p-1}^{a_r}, Y)$ doesn't have roots in $\ktl$.

This is not enough for our purposes, since we aim to prove that the specialization is irreducible. We can avoid the problem using proposition \ref{schinzel}; that is, we only need to prove that some auxiliary polynomial (call it $U$) doesn't have roots in the specialization. The problem is that $U$ isn't necessarily absolutely irreducible, so we need to work with each irreducible factor of $U$ at the same time. Suppose for simplicity that $r=0$ (it is not difficult to reduce to this case with a suitable change of variables), so $T$ and $U$ are polynomials in $X,Y$.

We can repeat the preceding construction to handle each irreducible factor of $U$, choosing the same prime $\mc{P}$ for each factor.

Let $U_0$ be an (absolutely) irreducible factor of $U$, and enlarge $k$ in order to ensure that $U_0 \in k[X,Y]$. Applying proposition \ref{lw} to $U_0$ we find and integer $a$ such that $U_0(\omega_{p-1}^{a}, Y)$ doesn't have roots in $\ktl$.  If we find the same integer $a$ for all irreducible factors, then $U(\omega_{p-1}^{a}, Y)$ itself doesn't have roots in $\ktl$, and finally $T(\omega_{p-1}^{a}, Y)$ is irreducible by proposition \ref{schinzel}. In general, though, our method will give a different values of $a$ for each factor $U_0$. This issue can be partially managed thanks to the following remark.

\begin{rem}
If $U_0(\omega_{p-1}^{a}, Y)$ does not have roots in $\ktl$, then the same will be true for each polynomial obtained by the action of $\Gal (\ktl/k)$. Since $U_0$ itself has coefficients in $k$, a conjugate will have the shape $U_0(\omega_{p-1}^{b a}, Y)$ for some suitable $b$ coprime with $p$. In our situation, knowing that $\Q^c \cap k \subset \Q(\omega_D)$, we can take every $b \equiv 1 \pmod{D}$.
\end{rem}

Thanks to this we are able to obtain the following lemma.

\begin{lemma} \label{coprimo}
Assume that each factor $U_0$ only depends on $X^D$. Moreover suppose that for each factor $U_0$ we can find some $a$ coprime with $m$ such that $U_0(\omega_{p-1}^{a}, Y)$ doesn't have roots in $\ktl$. Then $T(\omega_{p-1},Y)$ is irreducible over $\ktl$.
\end{lemma}

\begin{rem}
The assumption that $U_0$ depends only on $X^D$ may sound strange and quite restrictive at a first sight. Nevertheless we know that the polynomial we started with, namely $S_D$, has this property by construction. The problem lies in the fact that when one takes a factor of this polynomial, this property may be lost. Anyway one may hope to have some control, and for example to prove that $U_0$ only depends on $X^{D/D'}$, where $D'$ is little enough. The bigger is $D'$, the more delicate will be the estimates to carry out later.
\end{rem}

\begin{proof}[Proof of the lemma]
Fix a factor $U_0$ and consider the set $A \subset \Z/(p-1)\Z$ given by
\begin{equation*}
A = \left\{ a \in \Z/(p-1)\Z \text{ such that } U_0(\omega_{p-1}^{a}, Y) \text{ doesn't have roots in } \ktl \right\}.
\end{equation*}
Identify $\Z/(p-1)\Z$ with $\Z/m\Z \times \Z/D\Z$. Suppose that $A$ contains some $a$ coprime with $m$. Then by the remark $A \supset (\Z/m\Z)^{*} \times B$ for some $B \subset \Z/D\Z$. If moreover the polynomial $U_0$ only depends on $X^D$, then we can achieve $B = \Z/D\Z$.

In particular $(1,1) \in A$, so $U_0(\omega_{p-1}, Y)$ does not have roots in $\ktl$. Since this is true for each factor $U_0$, $U (\omega_{p-1}, Y)$ does not have roots in $\ktl$.
\end{proof}

At this point we face a problem: proposition \ref{lw} gives no control on whether $a$ is coprime with $m$, so we have to strengthen it a bit. Keep the identification
\begin{equation*}
\F_p^{*} \cong \Z/(p-1)\Z \cong \Z/m\Z \times \Z/D\Z.
\end{equation*}
The number of elements $x \in \F_p^{*}$ such that the projection on the first factor is coprime with $m$ is $D \varphi(m)$, where $\varphi$ is the Euler function. If we go back to the proof of proposition \ref{lw}, we want to compare this number with the upper bound in the estimate \eqref{stima}. Recall that we are dealing for simplicity with the case $r=0$, and that $\P$ has no inertia over $\Q$, so that $\mc{R}/\P=\F_p$, $p$ a prime. So we are able to obtain the stronger conclusion that $U_0(\omega_{p-1}^a, Y)$ is irreducible for some $a$ \emph{coprime with $m$} provided
\begin{equation*}
p \left( 1-\frac{d_{U_0} -1}{d_{\widetilde{U_0}}} \right) + O(p^{1/2}) \leq D \varphi(m).
\end{equation*}
Since $p-1=Dm$, what we need is an estimate from below for $\varphi (m)/m$. Remember that the existence of a prime $p$ with all the properties that we need is guaranteed by Chebotarev's theorem. Using an effective form of the theorem (such as in \cite{Lagarias}) one is able to bound $p$, and consequently $m$, from above. But we already know that $m$ doesn't have small prime factors, so this is translated in a bound for $\varphi (m)/m$.

Unfortunately this bound is not good enough for our purposes, but other tools from analytic number theory may do the trick. Once one is able to get this bound, the proof of proposition \ref{lw} shows that the hypothesis of lemma \ref{coprimo} can be fulfilled, and thus one gets a substantially different proof of the main arithmetical point in our proof.

\end{document}